\documentclass[12pt]{article}

\usepackage{amsmath}

\usepackage{amsthm}

\begin{document}

\makeatletter
\renewcommand\footnoterule{%
  \kern-3\p@
  \hrule\@width \textwidth
  \kern2.6\p@}
\renewcommand{\footnotesize}{\small}

\centerline{} 

\centerline{} 

\centerline {{\bf The Second Hankel Determinant for  }} 

\centerline{\bf  Starlike Functions of Order Alpha}

\centerline{} 

\centerline{\bf {D. K. Thomas}} 

\centerline{} 

\centerline{Department of Mathematics,} 
\centerline{Swansea University, Singleton Park,} 

\centerline{Swansea, SA2 8PP, UK.}

\centerline{d.k.thomas@swansea.ac.uk}

\renewcommand{\thefootnote}{}

\renewcommand{\baselinestretch}{1} 

\setcounter{page}{1} 

\setlength{\textheight}{21.6cm} 

\setlength{\textwidth}{14cm} 

\setlength{\oddsidemargin}{1cm} 

\setlength{\evensidemargin}{1cm} 

\pagestyle{myheadings} 

\thispagestyle{empty}

\markboth{\small{D K. Thomas}}{\small{}}

\date{}

\bigskip

\begin{abstract} Let $f$ be analytic in $D=\{z: |z|< 1\}$ with $f(z)=z+\sum_{n=2}^{\infty}a_{n}z^{n}$.  We give sharp bounds for  the second Hankel determinant  $H_{2}(2)=|a_{2}a_{4}-a_{3}^2|$ when $f$ is  starlike of order $\alpha$.

\end{abstract} 

\bigskip



\bigskip

\noindent 

{\bf Introduction, definitions and preliminaries} 

\bigskip

Let S be the class of analytic normalised univalent functions $f$, defined in $z\in D=\{z:|z|<1\}$ and given by

\begin{equation}
f(z)=z+\sum_{n=2}^{\infty}a_{n}z^{n}.
\end{equation}

\bigskip

The $qth$
Hankel determinant of $f$ is defined for $q\geq 1$ and $n\geq 0$ as follows, and has been extensively studied, see e.g. [1, 4, 5].
    \[
        H_{q}(n) = \left |
        \begin{array}{cccc}
        a_{n}\; & \mbox{$a_{n+1}...$}\;& \mbox{$a_{n+q+1}$}\\
        a_{n+1}\; &...& \mbox{$\vdots$}\\
        \vdots \\
        a_{n+q-1}\; & \mbox{$...$}\;& \mbox{$a_{n+2q-2}$}\\
        \end{array}
        \right |.
    \]
\bigskip

\bigskip

Denote by  $S^{*}$ the subclass of $S$ of  starlike functions, so that  $f\in S^{*}$ if, and only if, for $z\in D$

\begin{equation*}
Re  \ {\dfrac{zf'(z)}{f(z)}>0.}
\end{equation*}

\bigskip

Suppose that $f$ is analytic  in $D$  and given by $(1)$. Then $f$ is  starlike of order $\alpha$ in $D$ if, and only if, for $0\le \alpha <1$,
\bigskip

\begin{equation}
Re  \ {\dfrac{zf'(z)}{f(z)}>\alpha.}
\end{equation}

\bigskip

\noindent We denote this class by $S(\alpha)$,  so that $S(0)=S^*\subset S$ and $S^*(\alpha)\subset S^*$.

\bigskip

Let $P$ be the class of functions $p$ satisfying $Re \ p(z)>0$ for $z\in D$,  with $p(0)=1$.

\bigskip

\noindent Write

\begin{equation}
p(z)=1+\sum _{n=1}^{\infty}p_{n}z^{n}.
\end{equation}

\bigskip

We shall need the following result [3], which has been used widely.

\bigskip

\noindent {\bf Lemma}

\bigskip

Let $p\in P$ and be given by $(3)$, then for some complex valued $y$ with $|y|\le 1$, and some complex valued $\zeta$ with $|\zeta|\le 1$

\begin{equation*}
\begin{split}
2p_{2}&=p_{1}^{2}+y(4-p_{1}^{2}),\\
4p_{3}&=p_{1}^{3}+2(4-p_{1}^{2})p_{1}y-p_{1}(4-p_{1}^{2})y^{2}+2(4-p_{1}^{2})(1-|y|^{2})\zeta.\\
\end{split}
\end{equation*}
\bigskip

{\bf Results} 
\bigskip

It was shown in [2] that if $f\in S^*$ then $H_{2}(2)\le 1$,  and in [6] that if $f\in S(\alpha)$ with $0\le \alpha \le \dfrac{1}{2}$, then $H_{2}(2)\le (1-\alpha)^{2}$, both inequalities are sharp. 

\bigskip

\noindent We give the complete solution  for $f\in S(\alpha)$ when $0\le \alpha <1$ as follows.

\bigskip

\noindent {\bf Theorem}

\bigskip

\bigskip
If $f\in S^*(\alpha)$, then for $0\le \alpha <1$, the second Hankel determinant 

\bigskip

$$
H_{2}(2) \le (1-\alpha)^{2}.
$$

\bigskip

The inequality is sharp.

\bigskip

\begin{proof}

\bigskip

\noindent It follows from $(2)$ that we can write $zf'(z)=\alpha +(1-\alpha )f(z)p(z)$, and so equating coefficients we obtain

\bigskip

\begin{gather*}
a_{2}=(1-\alpha) p_{1}\\
a_{3}=\dfrac{1}{4} (2 (1 - \alpha)^2  p_{1}^2 + 2 p_{2} - 2 \alpha p_{2})\\
a_{4}=\dfrac{1}{6} (1 - \alpha) ((1 - \alpha)^2 p_{1}^3 +3 (1 - \alpha) p_{1} p_{2} + 2 p_{3}).\\
\end{gather*}

Hence

\bigskip

\begin{equation*}
\begin{split}
|a_{2}a_{4}-a_{3}^2|&=|-\dfrac{1}{48} (1 - \alpha )^2 (3 - 8 \alpha  + 4 \alpha ^2) p_{1}^4\\
 & - \dfrac{1}{4}p_{2}^2 + \dfrac{1}{2} \alpha p_{2}^2 - \dfrac{1}{4} \alpha^2p_{2}^2 + \dfrac{1}{3}p_{1} p_{3} - \dfrac{2}{3}
  \alpha p_{1} p_{3} + \dfrac{1}{3}  \alpha^2 p_{1} p_{3}|.
\end{split}
\end{equation*}

\bigskip

Noting that without loss in generality we can write $p_{1}=p$, with $0\le p\le 2$, we now use the Lemma to express the above in terms of $p$ to obtain   

\begin{equation*}
\begin{split}
|a_{2}a_{4}-a_{3}^2|&=|-\dfrac{1}{48} (1 - \alpha )^2 (3 - 8 \alpha  + 4 \alpha ^2) p^4\\
 & +\dfrac{1}{24} (1 - \alpha)^2 p^2 (4 - p^2) y - \dfrac{1}{12} (1 - \alpha)^2 p^2 (4 - p^2) y^2\\
 & -  \dfrac{1}{16} (1 - \alpha)^2 (4 - p^2)^2 y^2 + \dfrac{1}{6} (1 - \alpha)^2 p (4 - p^2) (1 - |y|^2)\zeta
|\\
&\le \dfrac{1}{48} (1 - \alpha )^2 |(3 - 8 \alpha  + 4 \alpha ^2)| p^4\\
  &+\dfrac{1}{24} (1 - \alpha)^2 p^2 (4 - p^2) |y| + \dfrac{1}{12} (1 - \alpha)^2 p^2 (4 - p^2) |y|^2\\
 & + \dfrac{1}{16} (1 - \alpha)^2 (4 - p^2)^2 |y|^2 + \dfrac{1}{6} (1 - \alpha)^2 p (4 - p^2) (1 - |y|^2):=\phi(|y|).
\end{split}
\end{equation*}
\bigskip

\bigskip

\noindent It is a simple exercise to show that $\phi'(|y|)\ge 0$  on [0, 1], so $\phi(|y|)\le \phi(1)$. Putting $|y|=1$ gives

\bigskip

\begin{equation*}
\begin{split}
|a_{2}a_{4}-a_{3}^2| &\le \dfrac{1}{48} (1 - \alpha )^2 |(3 - 8 \alpha  + 4 \alpha ^2)| p^4\\
  &+\dfrac{1}{8} (1 - \alpha)^2 p^2 (4 - p^2) + \dfrac{1}{16} (1 - \alpha)^2 (4 - p^2)^2\\
  &=1 - 2 \alpha + \alpha^2 - \dfrac{1}{16} (1 - \alpha)^2 p^4 + 
 \dfrac{1}{48} (1 - \alpha)^2 p^4 |3 - 8 \alpha + 4 \alpha^2|.
\end{split}
\end{equation*}

\bigskip

Considering $3 - 8 \alpha + 4 \alpha^2\ge 0$ and $3 - 8 \alpha + 4 \alpha^2\le 0$ separately, elementary calculus shows that the above expression is bounded by $(1-\alpha)^2$ in both cases.

\bigskip

We note that the  inequality in the Theorem is sharp when  $p_{1}=p_{3}=0$ and $p_{2}=2$.

\end{proof}

\end{document}